\documentclass{amsart}
\usepackage{amssymb}
\usepackage{amsmath}
\usepackage{amsthm}
\usepackage{mathptmx}
\usepackage{color}

\newcommand{\R}{\mathbb R}

\newcommand{\Z}{\mathbb Z}
\newcommand{\T}{\mathbb T}

\newtheorem{theorem}{Theorem}[section]

\newtheorem{corollary}[theorem]{Corollary}

\newtheorem*{remarks}{Remarks}
\newtheorem*{remarkwithout}{Remark}

\numberwithin{equation}{section}

\begin{document}
\title[Decay fifth order dispersive equations]{Decay properties for solutions of fifth order nonlinear dispersive equations}
\author{Pedro Isaza}
\address[P. Isaza]{Departamento  de Matem\'aticas\\
Universidad Nacional de Colombia\\ A. A. 3840, Medellin\\Colombia}
\email{pisaza@unal.edu.co}
\author{Felipe Linares}
\address[F. Linares]{IMPA\\
Instituto Matem\'atica Pura e Aplicada\\
Estrada Dona Castorina 110\\
22460-320, Rio de Janeiro, RJ\\Brazil}
\email{linares@impa.br}

\author{Gustavo Ponce}
\address[G. Ponce]{Department  of Mathematics\\
University of California\\
Santa Barbara, CA 93106\\
USA.}
\email{ponce@math.ucsb.edu}

\begin{abstract}
 We consider the initial value problem associated to a large class of fifth order nonlinear dispersive equations. This class includes several models arising in the study of different physical phenomena. Our aim is to establish special (space) decay properties of solutions to these systems. These properties complement previous unique continuation results and in some case, show that they are optimal. These decay estimates reflect the \lq\lq parabolic character" of these dispersive models in exponential weighted spaces. This principle was first obtained by T. Kato in solutions of the KdV equation.
\end{abstract}

\maketitle
\section{Introduction}

In this work we shall study decay and uniqueness properties of solutions to  a class of higher order dispersive models. More precisely,
we shall be concerned with one space dimensional (1D) dispersive models in which the dispersive relation is described by the fifth
order operator $\,\partial_x^5$. Roughly, the general form of the class of equations to be considered here is
\begin{equation}\label{I1}
\partial_t u-\partial_x^5u+P(u, \partial_x u, \partial_x^2u, \partial_x^3 u)=0
\end{equation}
where $P(\cdot)$ is a polynomial without constant or linear term, i.e.
\begin{equation}\label{I2}
P=P(x_1, x_2, x_3, x_4)=\underset{2\le |\alpha|\le N}{\sum}
a_{\alpha}\, x^{\alpha}, \hskip10pt N\in \Z^{+}, \hskip10pt N\ge 2,\hskip10pt a_{\alpha}\in\R.
\end{equation}

In this class one finds a large set of models arising in both mathematical and physical settings. Thus, the case
\begin{equation}\label{I2b}
P(u,\partial_xu, \partial_x^2 u, \partial_x^3 u)= 10 u\partial_x^3u+20 \partial_xu\partial_x^2u-30 u^2\partial_xu
\end{equation}
corresponds to the third equation in the KdV  hierarchy, where
\begin{equation}\label{I3}
\partial_t u +\partial_x u=0
\end{equation}
and the KdV
\begin{equation}\label{I4}
\partial_tu+\partial_x^3u+u\partial_xu=0
\end{equation}
are the first and second ones respectively in the hierarchy the jth being
\begin{equation}\label{I4b}
\partial_t u+ (-1)^{j} \partial_x^{2j+1}u +Q_j(u,\dots,\partial_x^{2j-1}u)=0,\hskip10pt j\in\Z^{+},
\end{equation}
with $Q_j(\cdot)$ an appropriate polynomial (see \cite{GGKM}).

Further examples of integrable models of the equation in \eqref{I1}-\eqref{I2} were deduced in \cite{kaup} and \cite{sawada-kotera} which also arise in the study of
higher order models of water waves.

In \cite{benney} the cases
\begin{equation}\label{I5}
P_1=c_1\,u\partial_xu, \hskip10pt P_2= u\partial_x^3u+ 2\partial_xu \partial_x^2u
\end{equation}
were proposed as models describing the interaction between long and short waves.

In \cite{lisher} the example of \eqref{I1}-\eqref{I2} with
\begin{equation}\label{I6}
P=(u+u^2)\partial_x u+(1+u)(\partial_xu\partial_x^2u+u\partial_x^3u)
\end{equation}
was deduced in the study of the motion of a lattice of anharmonic oscillators (by simplicity the values of the coefficients in \eqref{I6} have been taken equal to one).
In \cite{kakutani-ono} the equations \eqref{I1}-\eqref{I2} with $P=P_1$  as in \eqref{I5}  were proposed as a  model for magneto-acoustic waves at the critical angle in cold plasma.

Other cases of the equations \eqref{I1}-\eqref{I2} have been studied in \cite{olver}, \cite{kawahara}, \dots.

The well-posedness of the initial value problem (IVP) and the periodic  boundary value problem (PBVP) associated to the equation \eqref{I1} have been extensively studied (in the well-posedness
of these problems in a function space $X$ one includes existence, uniqueness of a solution $u\in C([0,T]:X)\cap \dots$ with $T=T(\|u_0\|_X)>0$, $u(\cdot,0)=u_0$, and the map $u_0\mapsto u$ being
locally continuous).

In \cite{saut} Saut proved the existence of solutions corresponding to smooth data for the IVP for the whole KdV hierarchy sequence of equations in \eqref{I4b}.

In \cite{schwarz} Schwarz considered the PBVP for the KdV hierarchy
\eqref{I4b} establishing existence and uniqueness in $H^s(\T)$ for
$s\ge 3j-1$.

In \cite{ponce}  Ponce showed that the IVP for \eqref{I1}-\eqref{I2} with \eqref{I2} as in \eqref{I2b} (i.e. the third equation in the KdV hierarchy) is globally well posed
in $H^s(\R)$ for $s\ge 4$.

In \cite{KPV1} Kenig, Ponce and Vega established the local well-posedness in weighted Sobolev spaces $H^s(\R)\cap L^2(|x|^m)$ of the IVP for the sequence of
equations in \eqref{I4b} for any polynomial $ Q_j(u, \dots, \partial_x^{2j}u)$, for $s\ge jm\,$ and $s\ge s_0(j)$.

The latter work  motivated several further studies concerning  the minimal regularity required on the Sobolev exponent $s$ such that the IVP \eqref{I1}-\eqref{I2} is
locally well-posed in $H^s(\R)$. These results heavily depend on the structure of the nonlinearity $P(\cdot)$ considered. For the precise statements of these kind of
results see \cite{kwon}, \cite{chen-guo}, \cite{kato}, \cite{chen-li-miao-wu}, \cite{huo} and references therein.  In particular, Kenig and Pilod \cite{kenig-pilod}
established local and global results in the energy space for the IVP associated to the equation \eqref{I1} with $P(\cdot)$ as in \eqref{I2b} (see also \cite{GKK}).

Special uniqueness properties of solutions to the IVP associated to the equation \eqref{I1}-\eqref{I2} were studied by  Dawson \cite{dawson}.

It was established in \cite{dawson} that if $u_1, u_2\in C([0,T]: H^6(\R)\cap L^2(|x|^3\,dx))$, $T>1$, are two solutions of  \eqref{I1}-\eqref{I2} such that
\begin{equation}\label{I8}
(u_1-u_2)(\cdot,0),\;(u_1-u_2)(\cdot,1)\in L^2\big(e^{x_{+}^{4/3+\epsilon}}\,dx\big)
\end{equation}
for some $\epsilon>0$, then $u_1\equiv u_2$. Moreover, in the case where in \eqref{I2} one has
\begin{equation*}
P=P(u,\partial_x u)=\underset{2\le \alpha_1+\alpha_2\le N}{\sum} a_{\alpha_1,\alpha_2}\, u^{\alpha_1} (\partial_xu)^{\alpha_2} \, , \hskip20pt N\in\Z^{+}
\end{equation*}
the exponent $4/3$ can be replaced by $5/4$.

In fact one should expect the general result in  \cite{dawson} to hold with $5/4$  in \eqref{I8} instead of $4/3$ for all $P(\cdot)$ in \eqref{I2}. However the argument of proof in
\cite{dawson} follows that given in \cite{EKPV} for the KdV equation. More precisely, it was established in \cite{EKPV} that  there exists $a_0>0$
such that if  $u_1, u_2\in C([0,T]:H^4(\R)\cap L^2(|x|^2\,dx))$, $T>1$, are solutions  of the IVP associated to the KdV equation \eqref{I4} with
\begin{equation}\label{I9}
(u_1-u_2)(\cdot,0),\;(u_1-u_2)(\cdot,1)\in L^2\big(e^{a_0\,x_{+}^{3/2}}\,dx\big)
\end{equation}
then $u_1\equiv u_2$. (Although, the statements in \cite{EKPV} and \cite{dawson} contain stronger hypotheses than the ones described in \eqref{I8} and \eqref{I9} respectively, these can be deduced by interpolation between \eqref{I8} and \eqref{I9} and the corresponding inequalities following by the assumptions on the class of solutions considered).

The value of the exponents above are dictated by the following decay estimate concerning the fundamental  solution of the associated linear problem
\begin{equation}\label{I10}
\begin{cases}
\partial_t v+ \partial_x^{2j+1}v=0,\\
v(x,0)=v_0(x).
\end{cases}
\end{equation}

In \cite{SSS} it was shown that
\begin{equation}\label{I11}
v(x,t)=\frac{c_j}{t^{1/(2j+1)}}\int\limits_{-\infty}^{\infty} K_j\Big(\frac{x-x'}{t^{1/(2j+1)}}\Big) v_0(x')\,dx'
\end{equation}
with $K_j(\cdot)$ satisfying
\begin{equation}\label{I12}
|K_j(x)|\le \frac{c}{(1+x_{-})^{(2j-1)/4j}}\, e^{-c_j\,x_{+}^{(2j+1)/2j}}
\end{equation}
with $x_{+}=\max \{x, 0\}$, $x_{-}=-\min\{x,\;0\}.$ Thus
\begin{equation}\label{I13a}
 K_j\Big(\frac{x}{t^{1/(2j+1)}}\Big)\sim e^{-c_j\,\big(\frac{x^{(2j+1)}}{t}\big)^{1/2j}}.
 \end{equation}

For this reason and the result in \cite{EKPV} one should expect that \eqref{I8} holds with $5/4$ instead of $4/3$ for a large class of polynomials in \eqref{I2}
including that in \eqref{I2b}. The obstruction appears in \cite{dawson} when the Carleman estimate deduced in \cite{EKPV} is used in this higher order setting.

In \cite{ILP} we proved that the result in \cite{EKPV} commented above is optimal. More precisely, the following results were established in \cite{ILP}.
\begin{itemize}
\item[(i)] If $u_0\in L^2(\R)\cap L^2(e^{a_0x_{+}^{3/2}}\,dx)$, $a_0>0$, then, for any $T>0$,  the solution of the IVP for the KdV equation \eqref{I4} satisfies
\begin{equation}\label{I13}
\underset{t\in[0,T]}{\sup} \int\limits_{-\infty}^{\infty} e^{a(t)\,x_{+}^{3/2}}\, |u(x,t)|^2\,dx \le  c^{*}=c^{*}(a_0;\,\|u_0\|_2;\,\|e^{\frac12a_0x_{+}^{3/2}}u_0\|_2; T),
\end{equation}
with
\begin{equation}\label{I14}
a(t)=\frac{a_0}{(1+27a_0^2\,t/4)^{1/2}}.
\end{equation}
\item[(ii)] If $u_1, u_2\in C([0,\infty) : H^1(\R)\cap L^2(|x|\,dx))$ are solutions of the IVP for the KdV equation \eqref{I4} such that
\begin{equation}\label{I15}
\int\limits_{-\infty}^{\infty} e^{a_0x_{+}^{3/2}} |u_1(x,0)-u_2(x,0)|^2\,dx<\infty,
\end{equation}
then, for any $T>0$, 
\begin{equation}\label{I16}
\underset{[0,T]}{\sup} \int\limits_{-\infty}^{\infty} e^{a(t)\,x_{+}^{3/2}} |u_1(x,t)-u_2(x,t)|^2\,dx\le c^{*}
\end{equation}
with
\begin{equation*}
\begin{split}
c^{*}= c^{*}(&a_0; \,\|u_1(\cdot,0)\|_{1,2};\,\|u_2(\cdot,0)\|_{1,2};\,\||x|^{1/2}u_1(x,0)\|_2; \\
&\;\;\|(u_1-u_2)(\cdot,0)\|_{1,2};  \|e^{\frac12a_0x_{+}^{3/2}}(u_1-u_2)(\cdot,0)\|_2; T)
\end{split}
\end{equation*}
with  $a(t)$ as in \eqref{I14}.
\end{itemize}

In order to simplify the exposition we state our main result first for the case of the IVP
\begin{equation}\label{I17}
\begin{cases}
\partial_tu-\partial_x^5u+b_1 u\,\partial_x^3u+ b_2 \partial_xu\,\partial_x^2u+ b_3u^2 \partial_xu=0,\\
u(x,0)=u_0(x),
\end{cases}
\end{equation}
with $b_1, \,b_2,\, b_3\in\R$ arbitrary constants.

\begin{theorem}\label{theorem1}
Let $a_0$ be a positive constant. For any given data
\begin{equation}\label{I18}
u_0\in H^4(\R)\cap L^2 (e^{a_0x_{+}^{5/4}}\,dx).
\end{equation}

The unique solution $u(\cdot)$ of the IVP  \eqref{I17}
\begin{equation*}
u\in C([0,T]: H^4(\R)) \cap \dots
\end{equation*}
satisfies
\begin{equation}\label{I19}
\underset{[0,T]}{\sup}\int\limits_{-\infty}^{\infty} e^{a(t) x_{+}^{5/4}} |u(x,t)|^2\,dx\le c^{*}=c^{*}(a_0; \|u_0\|_{4,2};\|e^{a_0 x_{+}^{5/4}}u_0\|_2; T)
\end{equation}
with
\begin{equation}\label{I20}
a(t)=\frac{a_0}{\sqrt[4]{1+k\,a_0^4t}} \text{\hskip15pt with \hskip10pt} k=11\,\frac{5^5}{4^5}.
\end{equation}
\end{theorem}


\begin{theorem}\label{theorem2}
Let $a_0$ be a positive constant. Let $u_1, u_2$ be solutions of the IVP \eqref{I17} such that
\begin{equation*}
\begin{split}
& u_1\in C([0,T]:H^8(\R)\cap L^2(|x|^4\,dx)),\\
&u_2\in C([0,T]: H^4(\R)) \cap L^2(|x|^2\,dx)).
\end{split}
\end{equation*}

If
\begin{equation}\label{I21}
\Lambda\equiv\int\limits_{-\infty}^{\infty} e^{a_0 x_{+}^{5/4}} |u_1(x,0)-u_2(x,0)|^2\,dx= \int\limits_{-\infty}^{\infty} e^{a_0 x_{+}^{5/4}} |u_{01}(x)-u_{02}(x)|^2\,dx <\infty,
\end{equation}
then, for $0<\epsilon\ll 1$
\begin{equation}\label{I22}
\underset{[0,T]}{\sup}\int\limits_{-\infty}^{\infty} e^{a(t)\, x_{+}^{5/4}} |u_1(x,t)-u_2(x,t)|^2\,dx\le c^{**}
\end{equation}
where $c^{**}= c^{**} (a_0;\, \|u_{01}\|_{8,2};\|u_{02}\|_{4,2}; \|x^2u_{01}\|_2;\, \|x\,u_{02}\|_{2}; \,\Lambda;\epsilon; T)$
and
\begin{equation*}
a(t)=\frac{a_0}{\sqrt[4]{1+k\,a_0^4t}} \text{\hskip20pt with \hskip10pt} k=k(\epsilon)= \frac{5^5}{4^5}\Big(\frac32+\frac{25}{4(5-\epsilon)}\Big).
\end{equation*}
\end{theorem}

\begin{remarkwithout}
In the case when the local solutions extend to  global ones, for example for the case of the model described in \eqref{I2b} for which the solutions satisfy infinitely
many conservation laws, the result in Theorem \ref{theorem1} holds in any time interval $[0,T]$. Indeed, this follows by combining the result in \cite{kenig-pilod} and integration by parts.
\end{remarkwithout}


Our next results generalize those in Theorems \ref{theorem1} and \ref{theorem2} to the following class of polynomials: 
\begin{equation}\label{I22b}
P(u, \partial_x u, \partial_x^2 u, \partial_x^3u)= Q_0(u, \partial_x u, \partial_x^2 u)\,\partial_x^3 u+  Q_1(u, \partial_x u, \partial_x^2 u)
\end{equation}
with
\begin{equation*}
Q_0(x_1,x_2,x_3)=\underset{1\le |\alpha|\le N}{\sum} a_{\alpha} x^{\alpha}, \hskip20pt N\in\Z^{+},\;\;N\ge 1, \;\; a_{\alpha}\in\R,
\end{equation*}
and 
\begin{equation*}
Q_1(x_1,x_2,x_3)=\underset{2\le |\alpha|\le M}{\sum} b_{\alpha} x^{\alpha}, \hskip20pt M\in\Z^{+},\;\;M\ge 2, \;\; b_{\alpha}\in\R.
\end{equation*}

Notice that all the nonlinearities in the models previously discussed belong to this class. For further discussion on the form of the
polynomial  $P(\cdot)$ in \eqref{I1}-\eqref{I2} see remark (iv) after the statements of Theorem \ref{theorem3} and Theorem \ref{theorem4}.

\begin{theorem}\label{theorem3} Let $a_0$ be a positive constant. For any given data
\begin{equation}\label{I23}
u_0\in H^{10}(\R)\cap L^2(e^{a_0\,x_{+}^{5/4}}\,dx).
\end{equation}
The unique solution $u(\cdot)$ of the IVP associated to the equation \eqref{I1} with $P(\cdot)$ as in \eqref{I22b}
\begin{equation}\label{I24}
u\in C([0,T]: H^{10}(\R))\cap \dots
\end{equation}
satisfies
\begin{equation}\label{I25}
\underset{[0,T]}{\sup}\int\limits_{-\infty}^{\infty} e^{a(t) x_{+}^{5/4}} |u(x,t)|^2\,dx\le c^{*}=c^{*}(a_0; \|u_0\|_{10,2};  \|e^{a_0 x_{+}^{5/4}}u_0\|_2;T)
\end{equation}
with
\begin{equation}\label{I26}
a(t)=\frac{a_0}{\sqrt[4]{1+k\,a_0^4t}} \text{\hskip15pt with \hskip10pt} k=11\,\frac{5^5}{4^5}.
\end{equation}
\end{theorem}

\begin{theorem}\label{theorem4} Let $a_0$ be a positive constant. Let $u_1, u_2$ be solutions of the IVP associated to the equation \eqref{I1} with $P(\cdot)$ as in \eqref{I22b}
in the class
\begin{equation*}
u_1, u_2 \in C([0,T]:H^{10}(\R)\cap L^2(\langle x\rangle^4\,dx))\cap \dots
\end{equation*}

If
\begin{equation}\label{I27}
\Lambda\equiv\int\limits_{-\infty}^{\infty} e^{a_0 x_{+}^{5/4}} |u_1(x,0)-u_2(x,0)|^2\,dx= \int\limits_{-\infty}^{\infty} e^{a_0 x_{+}^{5/4}} |u_{01}(x)-u_{02}(x)|^2\,dx <\infty,
\end{equation}
then for $0<\epsilon\ll 1$
\begin{equation}\label{I28}
\underset{[0,T]}{\sup}\int\limits_{-\infty}^{\infty} e^{a(t)\, x_{+}^{5/4}} |u_1(x,t)-u_2(x,t)|^2\,dx\le c^{**}
\end{equation}
where $c^{**}= c^{**} (a_0;\, \underset{j=1}{\overset{2}{\sum}}\big(\|u_{0j}\|_{10,2}+ \| x^2u_{0j}\|_2; \,\Lambda; \epsilon; T)$
and
\begin{equation*}
a(t)=\frac{a_0}{\sqrt[4]{1+k\,a_0^4t}} \text{\hskip20pt with \hskip10pt} k=k(\epsilon)= \frac{5^5}{4^5}\Big(\frac32+\frac{25}{4(5-\epsilon)}\Big).
\end{equation*}
\end{theorem}

\begin{remarks}\hskip20pt
\begin{itemize}
\item[(i)] The hypotheses on the regularity of the initial data (and consequently on the solutions) in Theorems \ref{theorem3} and \ref{theorem4}
are not sharp. They depend on the structure of the polynomial under consideration. However as written they apply to any polynomial $P(\cdot)$ in the class described in \eqref{I22b}.

\item[(ii)] Theorem \ref{theorem1} and Theorem \ref{theorem3} tell us that solutions of the equation \eqref{I1}  with $P(\cdot)$ as in \eqref{I2} decay accordingly to that
of the fundamental solution $K_2(\cdot)$ described in \eqref{I13a} with $j=2$.

\item[(iii)] Theorem \ref{theorem2} and Theorem \ref{theorem4} suggest that the results in \cite{dawson}  described in \eqref{I8} should hold with $5/4$ instead of $4/3$
in the exponent for any polynomial in \eqref{I2}. However this remains as an open problem.

\item[(iv)] We recall that in \cite{KPV1} a local existence theory in weighted Sobolev spaces ($H^s(\R)\cap L^2(\langle x\rangle^{2j}\,dx)$ with $s>4j$, $s$ large enough) for the
IVP associated to the equation \eqref{I1} with a general polynomial $P(\cdot)$ as in \eqref{I2} was established. This involves the use of a gauge transformation which transforms
the equation \eqref{I1} into an equivalent system. So one can ask if our argument presented here in Theorem \ref{theorem3} and Theorem \ref{theorem4} extends to this
general case. In this case, however, this general result  requires (decay and regularity) hypotheses involving the data (Theorem \ref{theorem3}) or the solutions (Theorem \ref{theorem4}) as well as some some of their derivatives. Also in this case the constant function $a(t)$ described in \eqref{I20} may be smaller (i.e. weaker decay).

\item[(v)] Concerning the existence of the solutions $u_1, u_2$ in the class described in Theorem \ref{theorem2} and Theorem \ref{theorem4} we recall the result in \cite{KPV1} . The fact that the operator $\Gamma=x+5t\partial_x^4$ commutes with $L= \partial_t -\partial_x^5$, and the \lq\lq identity''
\begin{equation*}
|x|^{\alpha}W(t)u_0= W(t) |x|^{\alpha}u_0 + W(t)\{\Phi_{t,\alpha}(\widehat u_0)(\xi)\}^{\lor}(x)
\end{equation*}
which holds for $\alpha\in(0,1)$,  where
 \begin{equation*}
\| \{ \Phi_{t,\alpha}(\widehat u_0)(\xi)\}^{\lor}\|_2\leq c(1+|t|)(\|u_0\|_2 + \| D^{4\alpha}u_0\|_2),
\end{equation*}
and  $W(t)$ denotes the unitary group associated to the linear equation in \eqref{I1} (see  \cite{FLP13}), imply in particular that in order to control $x^{\alpha}$ decay in the $L^2$ norm one needs at least to have $D^{4\alpha}_x$ derivatives in $L^2$. Thus
combining these ideas the result of existence of the solutions $u_1, u_2$ in the class described  follows.
\end{itemize}
\end{remarks}

In the proofs of Theorem \ref{theorem1} and Theorem \ref{theorem3} we need an intermediate decay result concerning the solutions of the IVP.
More precisely, in \cite{kato-83} T. Kato showed that $H^2$-solutions $u$ of the generalized KdV defined in the time interval $[0,T]$
\begin{equation*}
\partial_t u +\partial_x^3 u + u^k\,\partial_xu=0, \hskip20pt k\in\Z^{+},
\end{equation*}
with data $u_0\in L^2(e^{\beta x}\,dx)$, $\beta>0$ satisfy
\begin{equation}\label{I29}
e^{\beta x}\,u\in C([0,T]: L^2(e^{\beta x}\,dx))\cap C((0,T): H^{\infty}(\R)).
\end{equation}

Roughly in the linear case this follows from his observation that if $u$ is solution of 
\begin{equation*}
\partial_t u +\partial_x^3 u =0
\end{equation*}
and $v(x,t)= e^{\beta x} u(x,t)$, then
\begin{equation*}
\partial_t v+(\partial_x-\beta)^3 v =0.
\end{equation*}

So our next result, which will be used in the proofs of Theorems \ref{theorem1} and \ref{theorem3} extends Kato's result to solutions of the IVP \eqref{I1}
with $P(\cdot)$ as in \eqref{I22b}.

\begin{theorem}\label{theorem5}
Let $u\in C([0,T]: H^6(\R))$ be a solution of the IVP for the equation \eqref{I1} with $P$ as in \eqref{I22b}, corresponding
to data $u_0\in H^6(\R)\cap L^2(e^{\beta x}\,dx)$, $\beta>0$. Then
\begin{equation*}
e^{\beta x}u\in C([0,T]: L^2(\R))\cap C((0,T): H^{\infty}(\R)),
\end{equation*}
and
\begin{equation*}
\|e^{\beta x}u(t)\|_2\le c\,\|e^{\beta x}u_0\|_2, \;\;\; t\in [0,T].
\end{equation*}
\end{theorem}

We recall that although this is a subclass of the previously considered in Theorems \ref{theorem1}-\ref{theorem4}, this class includes all models
previously discussed. The restriction appears in our wish to use Kato's approach. In fact by using the idea developed in the proof of 
Theorems \ref{theorem1}-\ref{theorem4} one can extend the result in Theorem \ref{theorem5} to the whole class in \eqref{I2}.

Also the hypotheses in Theorem \ref{theorem5},  $u_0\in H^6(\R)$ can be significantly lower once  a particular form of the polynomial $P$ in \eqref{I22b}
is considered.

The paper is organized as follows. The construction of the weights to put forward the theory will be given in Section 1. The proofs of Theorem \ref{theorem1} and Theorem \ref{theorem2}
will be presented in Section 3. In Section 4, Theorem \ref{theorem3} and Theorem \ref{theorem4} will be proven. Finally, the proof of the extension of Kato's result will be detailed in Section 5.

\section{Construction of weights}
Consider the equation
\begin{equation}\label{A1}
\partial_t u-\partial_x^5u= F(x,t), \quad\quad t\ge 0,\;\;x\in\R.
\end{equation}

Formally, we perform (weighted) energy estimates in the equation \eqref{A1}, i.e., we multiply \eqref{A1} by 
$u\phi_N$, with $\phi_N=\phi_N(x,t)$ and $N\in\Z^{+}$, and integrate the result in the space variable. Thus after 
several integration by parts one gets
\begin{equation}\label{A2}
\begin{split}
\frac{d}{dt}\int u^2\phi_N\,dx& -\int u^2\,\partial_t\phi_N\,dx+5\int (\partial_x^2u)^2\,\partial_x\phi_N\,dx\\
&-5\int (\partial_xu)^2\,\partial_x^3\phi_N\,dx +\int u^2\partial_x^5\phi_N\,dx= 2\int Fu\phi_N\,dx.
\end{split}
\end{equation}

Using that
\begin{equation}\label{A3}
-5\int (\partial_xu)^2\partial_x^3\phi_N\,dx=5\int u\,\partial_x^2u\,\partial_x^3\phi_N\,dx-\frac52\int u^2\partial_x^5\phi_N \,dx.
\end{equation}

From the Cauchy-Schwarz and  Young inequalities  we have that for any $\epsilon \in [0,1]$
\begin{equation}\label{A4}
\begin{split}
5\Big| \int   u\,\partial_x^2u\,\partial_x^3\phi_N\,dx\Big |
 &\le 5\big(\int (\partial_x^2u)^2 \partial_x\phi_N\,dx\big)^{1/2}\big(\int u^2\frac{(\partial_x^3\phi_N)^2}{\partial_x\phi_N}\,dx\big)^{1/2}\\
 & \le (5-\epsilon)\int (\partial_x^2u)^2 \partial_x\phi_N\,dx + \frac{25}{4(5-\epsilon)}\int u^2\frac{(\partial_x^3\phi_N)^2}{\partial_x\phi_N}\,dx
 \end{split}
\end{equation}
(we remark that the integral above are taken on the set where $\partial_x^3\phi_N$ does not vanish. We will show that $\partial_x\phi_N$ does not vanish
in the support of  $\partial_x^3\phi_N$). Then, from \eqref{A2}-\eqref{A4} it follows that for any $\epsilon \in [0,1]$

\begin{equation}\label{A5}
\begin{split}
&\frac{d}{dt}\int u^2\phi_N\,dx -\int u^2\,\partial_t\phi_N\,dx+\epsilon\int (\partial_x^2u)^2\partial_x\phi_N\,dx\\
&-\frac32 \int u^2\partial_x^5\phi_N\,dx-\frac{25}{4(5-\epsilon)}\,\int u^2 \frac{(\partial_x^3\phi_N)^2}{\partial_x\phi_N}\,dx\\
&\le 2\int F u\phi_N\,dx,
\end{split}
\end{equation}
i.e. for $\epsilon\in[0,1]$
\begin{equation}\label{A6}
\begin{split}
&\frac{d}{dt}\int u^2\phi_N\,dx +\epsilon\int (\partial_x^2u)^2\partial_x\phi_N\,dx\\
&\le \int u^2\,\Big(\partial_t\phi_N+\frac32\partial_x^5\phi_N+\frac{25}{4(5-\epsilon)}\, \frac{(\partial_x^3\phi_N)^2}{\partial_x\phi_N}\Big)\,dx\\
&\;\;\;+2\int u F\phi_N\,dx.
\end{split}
\end{equation}

We shall use the inequality \eqref{A6} with $0<\epsilon\ll 1$. Then in order to simplify the proof we shall carry the details in the case $\epsilon=0$ and
remark that all the estimates involving the coefficient $25/4(5-\epsilon)$ are strict inequalities which also proves their extension to $\epsilon>0$ with $\epsilon\ll 1$.

We shall construct a sequence of weights $\{\phi_N\}_{N=1}^{\infty}$ which will be a key ingredient in the proof of our main theorems.

\begin{theorem}\label{theoremA}
Given $a_0>0$  and $\epsilon\in [0,1]$, $\epsilon\ll 1$, there exists a sequence $\{\phi_{\epsilon,N}\}_{N=1}^{\infty}\equiv\{\phi_N\}_{N=1}^{\infty}$ of functions with
\begin{equation}\label{A7}
\phi_N:\R\times[0,\infty)\to \R
\end{equation}
satisfying for any $N\in\Z^{+}$
\begin{enumerate}
\item[(i)] $\phi_N\in C^4(\R\times[0,\infty))$ \;\;with\;\; $\partial_x^5 \phi_N\phi(\cdot,t)$ having a jump discontinuity at $x=N$.
\item[(ii)] $\phi_N(x,t)>0$ \hskip10pt for all \;\;$(x,t)\in\R\times[0,\infty)$.
\item[(iii)] $\partial_x\phi_N(x,t)\ge 0$ \hskip10pt for all \;\; $(x,t)\in\R\times[0,\infty)$.
\item[(iv)]  There exist constants $c_N=c(N)>0$ and $c_0=c_0(a_0)>0$ such that
\begin{equation}\label{A8}
\phi_N(x,t)\le c_N\,c_0\,\langle x_{+} \rangle^4
\end{equation}
with
\begin{equation}\label{A9}
x_{+}=\max \{0;x\}, \hskip 15pt  \langle x \rangle=(1+x^2)^{1/2}.
\end{equation}
\item[(v)] For $T>0$ there is  $N_0\in\Z^{+}$ such that
\begin{equation}\label{A9b}
\phi_N(x,0) \le e^{a_0\,x_{+}^{5/4}}\text{\hskip10pt if \hskip10pt} N\ge N_0.
\end{equation}
Also
\begin{equation*}
\underset{N\uparrow\infty}{\lim} \phi_N(x,t) = e^{a(t)\,x_{+}^{5/4}}
\end{equation*}
for any $t>0$ and $x\in (-\infty,0)\cap(1,\infty)$ where
\begin{equation*}
a(t)=\frac{a_0}{\sqrt[4]{1+k\,a_0^4t}} \text{\hskip20pt with \hskip10pt} k=k(\epsilon)= \frac{5^5}{4^5}\Big(\frac32+\frac{25}{4(5-\epsilon)}\Big).
\end{equation*}
\item[(vi)] There exists a constant $c_0=c_0(a_0)>0$ such that for any $\epsilon\in[0,1]$, $\epsilon\ll 1$,
\begin{equation}\label{A10}
\partial_t\phi_N+\frac32\partial_x^5\phi_N+\frac{25}{4(5-\epsilon)}\, \frac{(\partial_x^3\phi_N)^2}{\partial_x\phi_N}\le c_0\,\phi_N
\end{equation}
for any $(x,t)\in\R\times[0,\infty)$.
\item[(vii)] There exist constants $c_j=c_j(j;a_0)>0$, $\;j=1,2,\dots,5$ such that
\begin{equation}\label{A11}
|\partial_x^j \phi_N(x,t)| \le c_j\,\langle x \rangle ^{j/4}\,\phi_N(x,t)
\end{equation}
for any $(x,t)\in\R\times[0,\infty)$.
\end{enumerate}
\end{theorem}

\vskip10pt

\noindent\underline{Proof of Theorem \eqref{theoremA}}
\vskip5pt

Given $a_0>0$,  for $N\in Z^{+}$ we define
\begin{equation}\label{A12}
\phi_N(x,t)=
\begin{cases}
e^{a(t)\varphi(x)}, \hskip20pt -\infty <x \le 1,\\
\\
e^{a(t)x^{5/4}}, \hskip20pt 1\le x\le N,\\
\\
P_N(x,t), \hskip20pt x\ge N,
\end{cases}
\end{equation}
where
\begin{equation}\label{A13}
a(t)=\frac{a_0}{\sqrt[4]{1+4k\,a_0^4\,t}} \le a_0, \hskip10pt t\ge 0,
\end{equation}
$a_0$ being the initial parameter and $k=k(\epsilon)>1$ is a constant whose precise value will be deduced below,
\begin{equation}\label{A14}
\varphi(x)= (1-\eta(x))\,x_{+}^5+\eta(x)\,x^{5/4},\;\;\;x_{+}=\max\{x;0\}
\end{equation}
for $x\in (-\infty,1]$ where $\eta\in C^{\infty}(\R)$, $\;\eta'\ge 0$ and
\begin{equation}\label{A15}
\eta(x)=
\begin{cases}
0, \hskip20pt x\le 1/2,\\
1, \hskip20pt x\ge 3/4,
\end{cases}
\end{equation}
(i.e. for each $x\in[0,1]$ $\varphi(x)$ is a convex combination of $x^5$ and $x^{5/4}$) and $P_N(x,t)$ is a polynomial of order $4$ in $x$ which matches the value of $e^{a(t)x^{5/4}}$
and its derivatives up to order $4$ at $x=N$:
\begin{equation}\label{A16}
\begin{split}
&P_N(x,t)=\\
&\hskip10pt \Big\{ 1+\frac54\, a N^{1/4}(x-N)+\frac{5}{4^2}\big(5 a^2 N^{2/4}+aN^{-3/4}\big)\frac{(x-N)^2}{2}\\
&\hskip15pt+ \frac{5}{4^3}\big( 25 a^3 N^{3/4}+15 a^2 N^{-2/4}- 3 a  N^{-7/4}\big)\frac{(x-N)^3}{3!}\\
&\hskip15pt+\frac{5}{4^4}\big( 125 a^4 N +150 a^3 N^{-1/4} -45 a^2N^{-6/4}+21 a N^{-11/4}\big)\frac{(x-N)^4}{4!}\Big\}\,e^{aN^{5/4}},
\end{split}
\end{equation}
with $a=a(t)$ as in \eqref{A13}.

Thus to prove \eqref{A8}-\eqref{A11} (i)-(vii) we consider the intervals $(-\infty,0]$, $[0,1]$, $[1, N]$ and $[N,\infty)$.

\vskip5mm
\noindent\underline{The interval $(-\infty,0]$}: In this case
\begin{equation*}
\phi_N(x,t)= e^{a(t)\cdot 0}= 1
\end{equation*}
which clearly satisfies  \eqref{A8} (i)-(vii).

\vskip5mm

\noindent\underline{The interval $[0,1]$}:  In this case
\begin{equation*}
\phi_N(x,t)= e^{a(t)\varphi(x)}
\end{equation*}
with
\begin{equation*}
\varphi(x)=(1-\eta(x)) \,x^5 +\eta(x)\,x^{5/4}\ge 0, \quad x\in [0,1]
\end{equation*}
with $\eta$ as in \eqref{A15}. Since in this interval $x^{5/4}\ge x^5$ it follows that
\begin{equation}\label{A17}
\begin{split}
\varphi'(x)&=(1-\eta(x))5x^4+\eta(x) \frac54 x^{1/4}+\eta'(x)(x^{5/4}-x^5)\\
&\ge (1-\eta(x))5x^4+\eta(x) \frac54 x^{1/4}\ge 0,
\end{split}
\end{equation}
and there exist $c_j>0$, $\; j=0,1,\dots,5$ such that
\begin{equation}\label{A18}
\varphi^{(j)}(x)\le c_j, \quad x\in[0,1].
\end{equation}
Since
\begin{equation}\label{A19}
a'(t)\le 0 \text{\hskip10pt one has \hskip10pt} a(t)\le a_0 \text{\hskip10pt for \hskip10pt} t\ge 0
\end{equation}
and we can conclude that
\begin{equation}\label{A20}
\partial_x^j\phi_N(x,t)\le c(j;a_0)\,\phi_N(x,t),\quad x\in [0,1], \;t\ge 0.
\end{equation}
Also
\begin{equation}\label{A21}
\partial_t\phi_N(x,t)=a'(t)\,\varphi(x)\,\phi_N(x,t) \le 0.
\end{equation}

Next we want to show that in this interval there exists $c_0=c_0(a_0)>0$ such that
\begin{equation}\label{A22}
\frac{(\partial_x^3\phi_N(x,t))^2}{\partial_x\phi_N(x,t)}\le c_0\,\phi_N(x,t),
\end{equation}
i.e.
\begin{equation}\label{A23}
(\partial_x^3\phi_N(x,t))^2\le c_0\,\phi_N(x,t)\, {\partial_x\phi_N(x,t)}.
\end{equation}

Since
\begin{equation*}
\begin{split}
\partial_x\phi_N&=a\varphi' \phi_N,\\
\partial_x^2\phi_N&=(a\varphi^{(2)}+(a\varphi')^2)\phi_N,\\
\partial_x^3\phi_N&=(a\varphi^{(3)}+3a^2 \varphi^{(2)} \varphi'+ (a\varphi')^3)\phi_N,
\end{split}
\end{equation*}
one has that for $x\sim 0$ \;$(x\ge 0)$
\begin{equation*}
\begin{split}
\partial_x\phi_N&\sim a5x^4 \phi_N,\\
\partial_x^2\phi_N&\sim (a 20 x^3+a^2 25 x^8)\phi_N,\\
\partial_x^3\phi_N&\sim (a 60 x^2 +3a^2 100 x^7+ a^3 125 x^{12})\phi_N,
\end{split}
\end{equation*}
Hence for $x\sim 0$ $\;(x\ge0)$
\begin{equation*}
\big(\partial_x^3\phi_N\big)^2\le c\,(a+a^3)^2\,x^4\,\phi_N^2
\end{equation*}
and
\begin{equation*}
\phi_N\,\partial_x\phi_N \ge 5 a \,x^4\,\phi_N^2.
\end{equation*}

Using \eqref{A19} (i.e. $a(t)\le a_0$ for $t\ge0$) it follows that there exists $\delta>0$ and a universal constant $c>0$ such that
\begin{equation}\label{A24}
\big(\partial_x^3\phi_N\big)^2\le c\,(a_0+a_0^5)\,\phi_N\,\partial_x\phi_N \quad \text{for\hskip10pt } x\in[0,\delta), \;t\ge 0.
\end{equation}

In the interval $[\delta,1]$ is easy to see that \eqref{A24} still holds (with a possible large $c>0$).

Combining the above estimates we see that \eqref{A8} (i)-(vii) hold in this interval.

\vskip5mm
\noindent\underline{The interval $[1,N]$}: In this region
\begin{equation}\label{A25}
\phi_N(x,t)=e^{a(t)\,x^{5/4}}, \hskip20pt x\in[1,N], \;\;t\ge 0.
\end{equation}

We calculate

\begin{equation}\label{A26}
\begin{split}
\partial_x\phi_N &= \frac54 \,ax^{1/4}\phi_N>0,\\
\partial_x^2\phi_N&=\frac5{4^2}\,(5a^2 x^{2/4}+a x^{-3/4})\phi_N,\\
\partial_x^3\phi_N&=\frac5{4^3}\,(25a^3x^{3/4}+15a^2x^{-2/4}-3ax^{-7/4})\phi_N,\\
\partial_x^4\phi_N&= \frac5{4^4}\,(125a^4x+150a^3x^{-1/4}-45a^2 x^{-6/4}+21a x^{-11/4})\,\phi_N,\\
\partial_x^5\phi_N&=\frac5{4^5}\,(625 a^5x^{5/4}+1250a^4-375 a^3 x^{-5/4}+375 a^2 x^{-10/4} -231 a x^{-15/4})\phi_N.
\end{split}
\end{equation}

Hence $\phi_N, \;\partial_x\phi_N >0$ and
\begin{equation}\label{A27}
\frac{\big(\partial_x^3\phi_N\big)^2}{\partial_x\phi_N}=\frac5{4^5}\Big( 625 a^5 x^{5/4}+750 a^4+75 a^3 x^{-5/4} - 90 a^2 x^{-10/4}+9 a x^{-15/4}\Big)\,\phi_N.
\end{equation}

Hence
\begin{equation}\label{A28}
\begin{split}
&\partial_t\phi_N+\frac32\partial_x^5\phi_N+\frac54\, \frac{\big(\partial_x^3\phi_N\big)^2}{\partial_x\phi_N}=\\
&\;\;=\Big\{ a'\,x^{5/4} +k\, a^5 x^{5/4}+c_4\,a^4+c_3\,a^3x^{-5/4}+c_2\,a^2 x^{-10/4}+c_1\,a x^{-15/4}\Big\}\,\phi_N
\end{split}
\end{equation}
with
\begin{equation}\label{A29}
\begin{split}
{\rm(a)} \hskip14pt k &= \frac{5^5}{4^5}\Big(\frac32+\frac54\Big)>1\\
{\rm(b)} \hskip10pt c_4 &= \frac{5}{4^5}\Big(\frac32 1250+\frac54 750\Big)>0\\
{\rm(c)} \hskip10pt c_3 &=\frac{5}{4^5} \Big(\frac32(-375)+\frac54 75\Big)< 0\\
{\rm(d)} \hskip10pt c_2 &=\frac{5}{4^5} \Big(\frac32 (375)+\frac54 (-90)\Big)>0\\
{\rm(e)} \hskip10pt c_1 &= \frac5{4^5}\Big(\frac32(-231)+\frac54\,9\Big)<0.
\end{split}
\end{equation}

Notice that if we change the coefficient $5/4$ in \eqref{A28} by $25/4(5-\epsilon)$, $\;\epsilon\in[0,1]$,  $\epsilon\ll 1$, the factor $5/4$ in \eqref{A29} (a)-(e) changes in a  similar
manner, and the value of $k$ in \eqref{A29} (a) will increase to
\begin{equation}\label{A30}
k(\epsilon)=\frac{5^5}{4^5}\Big(\frac32+\frac{25}{4(5-\epsilon)}\Big) >1
\end{equation}
and $c_1, c_2, c_3, c_4$ remain with the same sign, uniformly bounded in $\epsilon\in[0,1]$, $\epsilon\ll1$, and as we shall see below, the exact values
of $c_j$'s, $j=1,2,3,4$, are not relevant in the discussion below.

Next we solve the equation
\begin{equation}\label{A31}
a'(t)=- k\,a^5(t)
\end{equation}
which eliminates the terms with power $5/4$ on the right hand side of \eqref{A28}. Thus
\begin{equation}\label{A32}
a(t)=\frac{a_0}{\sqrt[4]{1+4\,k\,a_0^4\,t}}.
\end{equation}
Therefore to show that
\begin{equation}\label{A33}
\partial_t\phi_N+\frac32 \partial_x^5\phi_N+\frac54 \frac{\big(\partial_x^3\phi_N\big)^2}{\partial_x\phi_N}\le c_0\,\phi_N
\end{equation}
with $c_0=c_0(a_0)>0$  from \eqref{A29} it suffices to see that for $x\ge 1$.
\begin{equation}\label{A34}
c_4a^4+c_3 a^3 x^{-5/4}+c_2 a^2 x^{-10/4}+c_1 a x^{-15/4}\le c_0.
\end{equation}

Since $a(t)=a\le a_0$, $\;c_1, c_3\le 0$, and $x\ge 1$ one just needs to take $c_0$ such that
\begin{equation*}
c_4 a_0^4+c_2 a_0^2\le c_0.
\end{equation*}

Next, from \eqref{A26}
\begin{equation}\label{A35}
\begin{split}
\partial_x\phi_N &= \frac54 a x^{1/4}\phi_N\le c\,a_0\langle x\rangle^{1/4} \phi_N,\\
\partial_x^2\phi_N &\le c\,(a_0^2+a_0)\,\langle x\rangle^{1/2} \phi_N,\\
         \cdot \hskip10pt & \hskip40pt \cdot\\
          \cdot \hskip10pt & \hskip40pt \cdot\\
\partial_x^5\phi_N &\le c\,(a_0^5+a_0)\,\langle x\rangle^{5/4} \phi_N.
\end{split}
\end{equation}

Finally we remark that
\begin{equation}\label{A36}
\phi_N(x,t)=e^{a(t) x^{5/4} }\le e^{a_0\,N^{5/4}} \text{\hskip10pt for \hskip10pt} t\ge 0, \;\;x\in[1,N]
\end{equation}
which completes the proof of \eqref{A8} (i)-(vi) in this interval.

\vskip5mm
\noindent\underline{The interval $[N,\infty)$}: In this region
\begin{equation}\label{A37}
\begin{split}
&\phi_N(x,t)=\;P_N(x,t)= \Big\{ 1+\frac54 aN^{1/4}(x-N)\\
&+\frac{5}{4^2}\Big(\frac52 a^2 N^{2/4}+\frac12 aN^{-3/4}\Big) (x-N)^2\\
&+\frac{5}{4^3}\Big( \frac{25}6 a^3 N^{3/4}+\frac{15}{6} a^2N^{-2/4}-\frac{3}{6}a N^{-7/4}\Big)\, (x-N)^3\\
&+\frac{5}{4^4} \Big( \frac{125}{24} a^4 N^{4/4} +\frac{150}{24}a^3 N^{-1/4}-\frac{45}{24} a^2 N^{-6/4}+\frac{21}{24}a N^{-11/4}\Big)(x-N)^4\Big\}\,e^{aN^{5/4}},
\end{split}
\end{equation}
with $a=a(t)$ as in \eqref{A32}.

First we shall show that the negative coefficients of $(x-N)^3$, i.e. $-15aN^{-7/4}/384$ and of $(x-N)^4$, i.e. $-225a^2N^{-6/4}/6144$ can be controlled by the other ones. More precisely,
we shall see that there exists a universal constant $c>0$ such that
\begin{equation}\label{A34b}
\begin{split}
&\frac{5}{4^4} \Big( \frac{150}{24}a^3 N^{-1/4}-\frac{45}{24} a^2 N^{-6/4}+\frac{21}{24}a N^{-11/4}\Big)(x-N)^4\\
&-\frac5{4^3}\frac36 a N^{-7/4} (x-N)^3+\frac5{4^2}\frac12 a N^{-3/4} (x-N)^2\equiv R_N(x,t)\\
&\ge c\Big\{ (a^3N^{-1/4}+a^2N^{-6/4}+aN^{-11/4}) (x-N)^4\\
&\;\;\;\;\;\;\;\;+ aN^{-7/4}(x-N)^3 + a N^{-3/4}(x-N)^2\Big\},
\end{split}
\end{equation}

\begin{equation}\label{A35b}
\begin{split}
\partial_x R_N(x,t) \ge&\; c\Big\{ (a^3N^{-1/4}+a^2 N^{-6/4}+aN^{-11/4}) (x-N)^3\\
&\;\;\;\;+ aN^{-7/4}(x-N)^2+ a N^{-3/4} (x-N)\Big\},
\end{split}
\end{equation}
and
\begin{equation}\label{A36b}
\begin{split}
\frac{\partial_t R_N(x,t)}{a'(t)}\ge&\; c\Big\{ (a^2N^{-1/4}+aN^{-6/4}+N^{-11/4})(x-N)^4\\
&\;\;\;\;+ N^{-7/4}(x-N)^3+N^{-3/4}(x-N)^2\Big\}.
\end{split}
\end{equation}

Once \eqref{A34b}-\eqref{A36b} have been established it follows that there exists $c>0$ such that for $x\ge N$
\begin{equation}\label{A37b}
\begin{split}
P_N(x,t)\ge &c\Big\{ 1+ aN^{1/4}(x-N)+(a^2N^{2/4}+aN^{-3/4})\frac{(x-N)^2}{2!}\\
&+ (a^3N^{3/4}+a^2N^{-2/4}+a N^{-7/4})\frac{(x-N)^3}{3!}\\
&+( a^4N^{4/4}+a^3N^{-1/4}+a^2N^{-6/4}+aN^{-11/4})\frac{(x-N)^4}{4!}\Big\}\,e^{aN^{5/4}}\\
&\ge c\,e^{aN^{5/4}}>0,
\end{split}
\end{equation}
(which proves \eqref{A8} (ii) in this interval).

From \eqref{A37}-\eqref{A35b} it can also be seen that
\begin{equation}\label{A40a-A38}
\begin{split}
\partial_x P_N(x,t)\ge& \Big\{ \frac{5}{4} a N^{1/4}+\frac{5^2}{4^2} a^2\,N^{2/4} (x-N)\\
&+\frac{5^3}{4^3} \frac{a^3N^{3/4}}{2}(x-N)^2+\frac{5^4}{4^4} \frac{a^4N^{4/4}}{6}(x-N)^3\Big\}\;e^{aN^{5/4}}\ge 0
\end{split}
\end{equation}
(which proves \eqref{A8} (iii) in this interval)
and
\begin{equation}\label{A39}
\partial_t P_N(x,t)=a'(t)\, S_N(x,t)\,e^{aN^{5/4}}+a'(t)\,N^{5/4}\,P_N(x,t),
\end{equation}
where
\begin{equation}\label{A40}
\begin{split}
S_N(x,t)\ge&\; c\Big\{ N^{1/4}(x-N) +(aN^{2/4}+N^{-3/4})\frac{(x-N)^2}{2!}\\
&+( a^2N^{3/4}+aN^{-2/4}+a N^{-7/4})\frac{(x-N)^3}{3!}\\
&+(a^3N^{4/4}+a^2N^{-1/4}+aN^{-6/4}+a N^{-11/4})\frac{(x-N)^4}{4!}\Big\}\ge 0.
\end{split}
\end{equation}


The proof of  \eqref{A34b}  and \eqref{A35b} are similar, so we restrict ourselves to present the details of that for \eqref{A34b}.

First we observe that for any $\alpha>0$
\begin{equation}\label{A41}
a^2N^{-6/4}= \alpha a^{3/2}N^{-1/8}\,\frac{a^{1/2}N^{-11/8}}{\alpha}\le \frac12\Big(\alpha^2 a^3 N^{-1/4}+\frac{aN^{-11/4}}{\alpha^2}\Big).
\end{equation}

Thus taking $\alpha=\sqrt{6}$ it follows that
\begin{equation}\label{A42}
5(45+1) a^2 N^{-6/4}\le 5(149\,a^3N^{-1/4}+ 4\,aN^{-11/4}).
\end{equation}

Hence
\begin{equation}\label{A43}
\begin{split}
\frac{5}{4^4 24}\big(&150 a^3N^{-1/4}-45\,a^2N^{-6/4}+21 a N^{-11/4}\big)\\
&\ge \frac{5}{4^4 24}\big( a^3N^{-1/4}+a^2N^{-6/4}+17aN^{-11/4}\big).
\end{split}
\end{equation}

This takes care of the term with coefficient $\;-45\cdot5 a^2 N^{-6/4}/(4^4\cdot24)$ in $(x-N)^4$, see \eqref{A37b}. To handle the term $\; -15aN^{-7/4}(x-N)^3/(4^3\cdot6)$
 we write
\begin{equation}\label{A44}
\begin{split}
N^{-7/4}(x-N)^3=&\;\alpha\,N^{-3/8}\,(x-N)\,\frac{N^{-11/8}(x-N)^2}{\alpha}\\
\le &\; \frac12\Big(\alpha^2N^{-3/4}(x-N)^2+\frac{N^{-11/4}(x-N)^4}{\alpha^2}\Big)
\end{split}
\end{equation}
for any $\alpha>0$, and so taking $\alpha=\sqrt{6}$ one has that
\begin{equation}\label{A45}
\begin{split}
\frac{15+1}{4^3\cdot 6}aN^{-7/4}(x-N)^3&=\frac{1}{24}aN^{-7/4}(x-N)^3\\
&\le \frac18 aN^{-3/4}(x-N)^2+\frac{80}{4^4\cdot 24} a N^{-11/4}(x-N)^4.
\end{split}
\end{equation}

Collecting the above estimates we obtain \eqref{A34b}.

Next we shall show that there exists $c_0=c_0(a_0)>0$ (independent of $N$) such that if $x\ge N$
\begin{equation}\label{A46}
\partial_t\phi_N+\frac32\partial_x^5\phi_N+\frac{25}{4(5-\epsilon)}\frac{\big(\partial_x^3\phi_N)^2}{\partial_x\phi_N}\le c_0\,\phi_N,
\end{equation}
which in this region reduces to
\begin{equation}\label{A47}
\partial_t P_N\,\partial_xP_N+\frac{25}{4(5-\epsilon)}\big(\partial_x^3P_N)^2\le c_0\,P_N\,\partial_x P_N,
\end{equation}
for any $\epsilon\in[0,1]$, $\epsilon\ll 1$. As we have done before we first consider the case $\epsilon=0$.

Thus we have
\begin{equation}\label{A48}
\begin{split}
&\frac54\big(\partial_x^3P_N)^2=\frac54\Big\{\Big(\frac{5^3}{4^3} a^3N^{3/4}+\frac{75}{4^3} a^2N^{-2/4}-\frac{15}{4^3}aN^{-7/4}\Big)\\
&+\Big(\frac{5^4}{4^4}a^4N^{4/4}+\frac{750}{4^4}a^3N^{-1/4}-\frac{225}{4^4}a^2N^{-6/4}+\frac{105}{4^4}aN^{-11/4}\Big)(x-N)\Big\}^2\,e^{2aN^{5/4}},
\end{split}
\end{equation}
since, $a'(t)<0$, by \eqref{A39}  and \eqref{A37b} 
\begin{equation}\label{A49}
\partial_t P_N\le a'(t) N^{5/4} P_N \le a'(t) N^{5/4} \,e^{aN^{5/4}}.
\end{equation}
Thus, by \eqref{A40a-A38}
\begin{equation}\label{A50}
\begin{split}
\partial_t P_N\partial_xP_N &\le a'(t)\,N^{5/4}\times\\
 &\quad \;\;\Big\{\frac54 aN^{1/4}+\frac{5^2}{4^2} a^2N^{2/4}(x-N)+\frac{5^3}{4^3} a^3\frac{N^{3/4}}{2}(x-N)^2\Big\}\,e^{2aN^{5/4}}.
\end{split}
\end{equation}

First we shall use $\partial_t P_N\,\partial_xP_N$ to control the terms in \eqref{A48} involving the highest power in $N$. (Notice that we only handle the positive terms in \eqref{A48}).

Thus using \eqref{A29}--\eqref{A31} it follows that
\begin{equation}\label{A51}
\begin{split}
&\frac54 \frac{5^6}{4^6} a^6 N^{6/4}+ a'(t)N^{5/4}\frac54 aN^{1/4}= aN^{6/4}\frac54\Big(\frac{5^6}{4^6} a^5+a'(t)\Big)\\
&= aN^{6/4}\frac54\Big(\frac{5^6}{4^6} a^5-ka^5\Big)= a^6N^{6/4}\Big(\frac54 \frac{5^6}{4^6} -\frac{5^6}{4^6} \Big(\frac32+\frac{25}{4\cdot5}\Big)\Big)<0
\end{split}
\end{equation}
Notice that the last inequality above still holds with $25/4(5-\epsilon)$, $\epsilon\in[0,1]$, $\epsilon\ll 1$, instead of $25/4\cdot5$ and $5/4$.

Also by \eqref{A29}--\eqref{A31}
\begin{equation}\label{A52}
\begin{split}
&\frac54 2 \Big(\frac{5^3}{4^3} a^3 N^{3/4}\Big)\Big(\frac{5^4}{4^4} a^4N^{4/4}\Big)(x-N)+a'(t)N^{5/4}\frac{5^2}{4^2}a^2N^{2/4}(x-N)\\
&=a^2N^{7/4}\frac{5^2}{4^2}(x-N)\Big(2\frac{5^6}{4^6}a^5+a'(t)\Big)\\
&=a^2N^{7/4}\frac{5^2}{4^2}(x-N)\Big(2\frac{5^6}{4^6}a^5-ka^5\Big)\\
&=a^7N^{7/4}\frac{5^2}{4^2}(x-N)\Big(2\frac{5^6}{4^6}-\frac{5^5}{4^5}\Big(\frac32+\frac{25}{20}\Big)\Big)\le0
\end{split}
\end{equation}
(where the remark after \eqref{A51} also applies), and again by \eqref{A29}--\eqref{A31}
\begin{equation}\label{A53}
\begin{split}
&\frac54\frac{5^8}{4^8}a^8N^2(x-N)^2+a'(t)N^{5/4}\frac{5^3}{4^3}\frac{a^3N^{3/4}}{2}(x-N)^2\\
&=\frac{5^3}{4^3}a^3N^2(x-N)^2\Big(\frac{5^6}{4^6}a^5+\frac12 a'(t)\Big)\\
&=\frac{5^3}{4^3}a^3N^2(x-N)^2\Big(\frac{5^6}{4^6}a^5-\frac12 ka^5\Big)\\
&=\frac{5^3}{4^3}a^8N^2(x-N)^2\Big(\frac{5^6}{4^6}-\frac12\frac{5^5}{4^5}\Big(\frac32+\frac{25}{20}\Big)\Big)\le0
\end{split}
\end{equation}
(where the remark after \eqref{A51} also applies).

We bound the remaining terms in \eqref{A48} by $c_0 \partial_tP_N\partial_xP_N$ . For that we use the fact that
\begin{equation*}
e^{aN^{5/4}}\partial_x P_N \le P_N\,\partial_xP_N, \quad x\ge N.
\end{equation*}

Thus, from \eqref{A40a-A38} 
\begin{equation}\label{A54}
\begin{split}
\frac54 2\frac{5^4}{4^4} a^4N^{4/4}\frac{750}{4^4} a^3N^{-1/4} (x-N)^2e^{2aN^{5/4}}
&\le c_0 \frac{5^3}{4^3} a^3 \frac{N^{3/4}}{2}(x-N)^2e^{2aN^{5/4}}\\
&\le c_0\, \partial_xP_N e^{aN^{5/4}}\le c_0\, P_N  \partial_xP_N ,
\end{split}
\end{equation}
by taking
\begin{equation}\label{A55}
c_0> ca_0^4,\quad  c \text{\hskip10pt universal constant}.
\end{equation}

Notice that with this choice of $c_0$ \eqref{A54} holds even when the factor $5/4$ in the left hand side is replaced by $25/4(5-\epsilon)$,
$\;\epsilon\in[0,1]$,$\;\epsilon\ll1$. Also, from \eqref{A40a-A38}
\begin{equation}\label{A56}
\begin{split}
&2\frac54\Big(\frac{5^4}{4^4}a^4N^{4/4}\frac{75}{4^3} a^2N^{-2/4}+\frac{750}{4^4} a^3 N^{-1/4}\frac{5^3}{4^3} a^3 N^{3/4}\Big)(x-N)e^{2aN^{5/4}}\\
&\le c_0\frac{5^2}{4^2} a^2N^{2/4}\,(x-N)\, e^{2aN^{5/4}}\le c_0\,\partial_xP_N\,e^{aN^{5/4}}\le c_0\,P_N \partial_x P_N
\end{split}
\end{equation}
by taking $c_0$ as in \eqref{A55} (where the remark after \eqref{A55} also applies). Also
\begin{equation}\label{A57}
2\frac54\frac{5^3}{4^3} a^3 N^{3/4} \frac{75}{4^3} a^2 N^{-2/4}e^{2aN^{5/4}} \le c_0 \frac{5}{4} a N^{1/4}e^{2aN^{5/4}}\le  c_0\,P_N \partial_x P_N
\end{equation}
by taking $c_0$ as in \eqref{A55} (and the remark after \eqref{A55} also applies).

This handles all the terms in \eqref{A48} having a positive coefficient and a positive power of $N$. The reminder ones having positive coefficients
can be bounded by
\begin{equation*}
c_0\,\frac54 a N^{1/4}.
\end{equation*}

Combining the above estimates with \eqref{A34b}--\eqref{A40} completes the proof of \eqref{A47}.

Finally \eqref{A37b} yields \eqref{A11} in this region $x\ge N$.

To finish the proof  we need to prove (v) in the region $[N,\infty)$. We use \eqref{A26}   with $t=0$   and observe that the negative terms in the expression for $\frac{d^5}{dx^5}e^{a_0x^{5/4}}$ can be absorved by the positive terms for $x\geq N$ and $N$ sufficiently large. More precisely, 
\begin{equation*}
1250a_0^4>2\cdot 375a_0^3x^{-5/4}
\quad\text{and }\;\;\;\;
{375}a_0^2N^{-10/4}>2\cdot 231 a_0 x^{-15/4}
\end{equation*} 
if  $x^{5/4}>c/a_0$,  where $c$ is an absolute constant. To have this for $x\geq N$, it sufficies to take $N$ in such a way that 
$N^{5/4}>c/a_0$. This is, $N>c^{4/5}a_0^{-4/5}\equiv N_0$.

 In this way, 
\begin{equation*}
\frac{d^5}{dx^5}(e^{a_0x^{5/4}}-P_N(x,0))= \frac{d^5}{dx^5}e^{a_0x^{5/4}}\geq 0
\end{equation*}
 for $x> N>N_0$. Since $e^{a_0x^{5/4}}$ and $P_N(x,0)$ coincide at $x=N$ up to the fourth derivative, we conclude that 
$e^{a_0x^{5/4}}\geq P_N(x,0)$ for $x\geq N\geq N_0$ , which proves (v) in this region.

Thus we have completed the proof of \eqref{A7} (i)-(vii), \eqref{A8}--\eqref{A11}.

\begin{corollary} \label{corollaryA} There exists $\widetilde{c}_0=\widetilde{c}_0(a_0; T)>0$  such that for any $N\in\Z^{+}$ sufficiently large, $x\in\R$, $t\in[0,T]$
\begin{equation}\label{A58}
\phi_N(x,t) \le  \widetilde{c}_0\,\big(1+\langle x \rangle \,\partial_x\phi_N(x,t)\big).
\end{equation}
\end{corollary}

The proof follows from the construction of the weight $\phi_N$.

\section{Proofs of Theorem \ref{theorem1} and Theorem \ref{theorem2}}

\noindent \underline{Proof of Theorem \ref{theorem1}}
\vskip5pt

Using the result in Theorem \ref{theorem5}
(and the remark afterwards) we have that our solution $u$ of the IVP \eqref{I17}
satisfies
\begin{equation}\label{3.2}
u\in C([0,T];H^4(\mathbb R)\cap L^2(e^{\beta x}dx))\quad\text{for any } \beta>0.
\end{equation}
Therefore, by interpolation one has that
\begin{equation}\label{3.2b}
\partial_x^ju\in C ([0,T];H^{4-j}(\mathbb R)\cap L^2(e^{(4-j)\beta x/4} dx))\,\quad j=0,1,2,3,4.
\end{equation}
In particular $u\in C([0,T];L^2(\langle x\rangle ^kdx))$, for any $k$. Suppose first that $u$ is sufficiently regular, say $u\in C([0,T];H^5(\mathbb R))$. Then we can perform energy estimates for $u$ using the weights $\{\phi_N\}$  (since $\phi_N\leq c\langle x \rangle ^4$).
Thus, we multiply the equation in \eqref{I17} by $u\phi_N$ and integrate by parts in the space variable to obtain 
\begin{equation}\label{3.1}
\int\partial_tuu\phi_N-\int\partial_x^5uu\phi_N+b_1\int u\partial_x^3uu\phi_N+b_2\int \partial_xu\partial_x^2uu\phi_N+b_3\int u^2\partial_xuu\phi_N=0,
\end{equation}
and applying \eqref{A6} and \eqref{A10} we have
\begin{equation}\label{3a}
\begin{aligned}
2\bigl(\int\partial_tuu\phi_N\,dx-&\int\partial_x^5uu\phi_N\,dx\bigr)\geq \frac{d}{dt}\int u^2\phi_N\,dx+\varepsilon\int(\partial_x^2u)^2\partial_x\phi_N\,dx\\
&-\int u^2\Bigl(\partial_t\phi_N+\frac32\partial_x^5\phi_N+\frac{25}{4(5-\varepsilon)}\frac{(\partial_x^3\phi_N)^2}{\partial_x\phi_N}\Bigr)\,dx\\
&\geq  \frac{d}{dt}\int u^2\phi_N\,dx+\varepsilon \int(\partial_x^2u)^2\partial_x\phi_N\,dx-c_0\int u^2\phi_N\,dx,
\end{aligned}
\end{equation}
with $\varepsilon\in[0,1]$, $\varepsilon<<1$, and $c_0=c_0(a_0)$. 
In  the proof of Theorem \ref{theorem1} we will only use \eqref{3a} with $\varepsilon=0$.

Now we shall handle the third, fourth and fifth terms on the right hand of  \eqref{3.1}. Thus we write
\begin{equation}\label{3.6}
\int u \partial_x^3u\, u\phi_N\,dx \le c\,\|\partial_x^3u\|_{\infty}\int u^2\phi_N\,dx,
\end{equation}
by integration by parts
\begin{equation}\label{3.7}
\begin{split}
\int \partial_x u\partial_x^2u\,u\phi_N\,dx&=-\frac12\int \partial_x^3u \,u^2\phi_N\,dx -\frac12\int \partial_x^2u u^2\partial_x\phi_N\,dx\\
&\equiv E_1+E_2.
\end{split}
\end{equation}

The bound for $E_1$ is similar to that in \eqref{3.6}. To control $E_2$ we recall that (see \eqref{A11})
\begin{equation}\label{3.8}
0\le \partial_x\phi_N(x,t)\le c_1\,\langle x\rangle^{1/4} \phi_N(x,t)\le c(1+e^x)\,\phi_N(x,t)
\end{equation}
so
\begin{equation*}
E_2\le c\,\big(\|e^x\,\partial_x^2u\|_{\infty}+\|\partial_x^2u\|_{\infty}\big)\,\int u^2\phi_N\,dx.
\end{equation*}

Notice that by combining Sobolev embedding and \eqref{3.2b} one has that $\int_0^T \| e^{x}\partial_x^2u\|_{\infty}(t)\,dt$ is finite.

Finally for the fifth term in \eqref{3.1} we have that
\begin{equation}\label{3.9}
\int u^2\partial_xu\,u\phi_N\,dx \le \|u\,\partial_x u\|_{\infty}\int u^2\phi_N\,dx.
\end{equation}

Collecting the above information, from \eqref{3a}  we can conclude that for any $N\in\Z^{+}$
\begin{equation*}
\frac{d}{dt} \int u^2(x,t)\phi_N(x,t)\,dx \le M(t)\,\int u^2(x,t)\phi_N(x,t)\,dx
\end{equation*}
with $M(t)\in L^{\infty}([0,T])$, where $M(\cdot)$ depends on $\, a_0, \,\| e^{x}\,u_0\|_2, \,\|u_0\|_{4,2}$. Hence, from property (v) in \eqref{A9b}, and Gronwall's Lemma we see that for $t\in [0,T]$
\begin{equation}\label{G}\begin{aligned}
\int u^2(x,t)\phi_N(x,t)\,dx&\leq c\bigl(\int u_0^2(x)\phi_N(x,0)\,dx\bigr)e^{\int_0^TM(t')\,dt'}\\
&\leq c(a_0, \|e^{\frac12a_0^2x_+^{5/4}}u_0\|_2, \| u_0\|_{4,2},T)\int u_0^2(x)e^{a_0x_+^{5/4}}\,dx\,.
\end{aligned}
\end{equation}

Now, we will establish \eqref{G} for our less regular solution $u\in C([0,T]; H^4(\mathbb R))$. To do that, we  consider the IVP \eqref{I17} with 
regularized initial data $u_{0,\delta}:=\rho_\delta*u(\cdot+\delta,0)$, where $\delta>0$,  $\rho_\delta=\frac1\delta\rho(\frac{\cdot}\delta)$, $\rho\in C^\infty(\mathbb R)$ is supported in (-1,1), and $\int\rho=1$. 
Since 
\begin{equation}\label{3c}
u_{0,\delta}\to u_0 \text{\hskip5pt  in \hskip5pt } H^4(\mathbb R)\text{\hskip5pt  as \hskip5pt } \delta\to 0,
\end{equation} 
by the well-posedness result in \cite{ponce} for  the IVP \eqref{I17} in $H^4(\mathbb R)$,  the corresponding solutions $u_\delta$ satisfy $u_\delta(t)\to u(t)$ in $H^4(\mathbb R)$ uniformly for $t\in [0,T]$ as $\delta\to 0$. In particular, by Sobolev embeddings, for fixed $t$ 
\begin{equation}u_{\delta}(x,t)\to u(x,t)\quad \text{for all }x\in\mathbb R\quad\text{as }\delta\to 0.\label{3d}\end{equation} 
Also, it can be proved (see Theorem 1.1 in \cite{ILP}) that
\begin{equation}
\|e^{\frac12a_0x_+^{5/4}}u_{0,\delta}\|_2\leq \|e^{\frac12a_0x_+^{5/4}}u_{0}\|_2.\label{3e}
\end{equation}
Since $u_\delta$ is sufficiently regular we have \eqref{G} with $u_\delta$ and $u_{0,\delta}$ instead of $u$ and $u_0$.
In this way, for $t$ fixed,   using  \eqref{3c}-\eqref{3e}, and applying Fatou's Lemma we see that
\begin{equation}
\int u^2(x,t)\phi_N(x,t)\,dx\leq C(a_0, \|e^{\frac12a_0x_+^{5/4}}u_0\|_2, \|u_0\|_{4,2},T)\int u_0^2(x)e^{a_0x_+^{5/4}}\,dx
\end{equation}
Now, we make $N\to \infty$ and apply  property (v) in Theorem \ref{theoremA} and Fatou's Lemma again to obtain

\begin{equation*}
\sup_{t\in[0,T]}\int u^2(x,t)e^{a(t)x_+^{5/4}}\,dx\leq c*,
\end{equation*}
which is the desired result.

\vskip10pt

\noindent\underline{Proof of Theorem \ref{theorem2}}
\vskip5pt

We consider the equation for the difference of the two solutions
\begin{equation}\label{3.10}
w(x,t)=(u_1-u_2)(x,t)
\end{equation}
that is,
\begin{equation}\label{3.11}
\begin{split}
\partial_t w-\partial_x^5 w=& \; -b_1(u_1\partial_x^3w+\partial_x^3u_2 w)-b_2(\partial_x u_1\partial_x^2w+\partial_x^2u_2\partial_xw)\\
&- b_3(\partial_xu_2(u_1+u_2)w+u_1^2\partial_x w).
\end{split}
\end{equation}

We follow the argument given in the proof of Theorem \ref{theorem1} with $\epsilon\in[0,1]$, $\epsilon\ll 1$. Hence we multiply \eqref{3.11} by $w\phi_N$ and integrate in the
variable $x$ and use that
\begin{equation}\label{3.12}
\int u_1 \partial_x^3 w \,w\phi_N\,dx= \frac12 \int u_1\phi_N \partial_x^3(w^2)\,dx-3\int u_1\phi_N \partial_x w\partial_x^2w \,dx\equiv F_1+F_2
\end{equation}
where
\begin{equation}\label{3.13}
F_1=-\frac12 \int \partial_x^3(u_1\phi_N) w^2\,dx.
\end{equation}
Then using \eqref{A11} it follows that
\begin{equation}\label{3.14}
|F_1|\le \underset{j=0}{\overset{3}{\sum}} \|\langle x\rangle ^{j/4}\partial_x^{3-j} u_1\|_{\infty} \int w^2 \phi_N\,dx
\end{equation}
and after some integration by parts
\begin{equation}\label{3.15}
\begin{split}
F_2=-\frac32 \int u_1\phi_N \partial_x(\partial_xw)^2\,dx
&=-\frac32 \int \partial_x(u_1\phi_N) w\partial_x^2w\,dx +\frac34\int \partial_x^3(u_1\phi_N)w^2\,dx\\
&\equiv F_2^1+F_2^2.
\end{split}
\end{equation}

We observe that the same bound for $F_1$ given in \eqref{3.14} applies to $F_2^2$. For $F_2^1$ we write
\begin{equation}\label{3.16}
F_2^1=-\frac32 \int \partial_x u_1\phi_N w\partial_x^2w\,dx-\frac32 \int u_1\partial_x\phi_N w \partial_x^2w\,dx\equiv F^{1,1}_2+F^{1,2}_2,
\end{equation}
with
\begin{equation}\label{3.17}
\begin{split}
|F^{1,2}_2| &\le \frac{\epsilon}4\int (\partial_x^2 w)^2\partial_x\phi_N\,dx+\frac{4}{\epsilon}\int u^2_1 w^2 \partial_x\phi_N\,dx\\
&\le \frac{\epsilon}4\int (\partial_x^2 w)^2\partial_x\phi_N\,dx+\frac{c}{\epsilon}\|u_1\langle x\rangle^{1/8}\|_{\infty}^2\int w^2\phi_N\,dx
\end{split}
\end{equation}
using \eqref{A11} and
\begin{equation*}
\begin{split}
|F^{1,1}_2| &\le c_0\int |\partial_x u_1(1+\langle x\rangle \partial_x\phi_N)w\partial_x^2 w|\,dx\\
&=c_0\,\int |\partial_x u_1 w\partial_x^2 w|\,dx+c_0\,\int |\partial_x u_1\langle x\rangle w\partial_x^2 w|  \partial_x\phi_N \,dx\\
&\le c_0\,\int |\partial_x u_1w\partial_x^2 w|\,dx+\frac{\epsilon}4\int (\partial_x^2w)^2\partial_x\phi_N\,dx
+\frac{c_0'}{\epsilon}\int |\partial_x u_1|^2 \langle x\rangle^2 w^2\partial_x\phi_N\,dx\\
&\le c_0\,\int  |\partial_x u_1w\partial_x^2 w|\,dx+\frac{\epsilon}4\int (\partial_x^2w)^2\partial_x\phi_N\,dx
+\frac{c}{\epsilon} \|\langle x\rangle^{1+1/8}\partial_x u_1\|_{\infty}^2\int  w^2\phi_N\,dx
\end{split}
\end{equation*}
by using Corollary \ref{corollaryA} \eqref{A58} and \eqref{A11}.

Directly one has that
\begin{equation*}
\int \partial_x^3 u_2 w^2 \phi_N\,dx \le \|\partial_x^3 u_2\|_{\infty} \int w^2\phi_N\,dx.
\end{equation*}

The estimate for the term
\begin{equation*}
\int \partial_x u_1\partial_x^2w w\phi_N\,dx
\end{equation*}
is similar to that given above for $F_2^{1,1}$.

Similarly, we have that
\begin{equation*}
\int \partial_x^2 u_2\partial_xw w \phi_N\,dx=-\frac12\int \partial_x(\partial_x^2u_2\phi_N)w^2\,dx
=-\frac12 \int \partial_x^3u_2 w^2 \phi_N\,dx -\frac12\int \partial_x^2u_2\partial_x\phi_N w^2\,dx
\end{equation*}
with
\begin{equation*}
\big|\int \partial_x^2u_2\partial_x\phi_N w^2\,dx\big| \le\int |\partial_x^2u_2|\langle x\rangle^{1/4}\phi_N w^2\,dx
\le \| \langle x\rangle^{1/4}\partial_x^2u_2\|_{\infty}\int w^2\phi_N\,dx
\end{equation*}
and 
\begin{equation*}
\big|  \int \partial_x^3u_2 w^2 \phi_N\,dx\big | \le \|\partial_x^3 u_2\|_{\infty} \int w^2\phi_N\,dx.
\end{equation*}

Finally, the terms
\begin{equation*}
\int \partial_xu_2(u_1+u_2)w^2\phi_n\,dx +\int u_1^2\partial_x w w\phi_N\,dx
\end{equation*}
can be handled analogously.

Thus combining the inequalities 
\begin{equation*}
\int \partial_t w w\phi_N\,dx -\int \partial_x^5w w\phi_N\,dx \ge \; 2 \frac{d}{dt}\int w^2\phi_N\,dx +\epsilon\int (\partial_x^2 w)^2\partial_x\phi_N\,dx
-c_0\int w^2\,\phi_N\,dx
\end{equation*}
(see \eqref{3a}),
\begin{equation*}
\int |\partial_x u_2 w\partial_x^2 w|\,dx \le \|\partial_x u_2\|_{\infty}\|w\|_2\|\partial_x^2 w\|_2\equiv L(t),
\end{equation*}
and the above estimates we have that
\begin{equation*}
\frac{d}{dt}\int w^2(x,t)\phi_N(x,t)\,dx \le M(t) \int  w^2(x,t)\phi_N(x,t)\,dx + L(t)
\end{equation*}
where 
\begin{equation*}
\begin{split}
M(t)=& \; c(\epsilon)\Big(\underset{j=0}{\overset{3}{\sum}}\|\langle x\rangle^{j/4}\partial_x^{3-j}u_1\|_{\infty}
+\|\langle x\rangle^{1+1/8}\partial_x u_1\|_{\infty}\\
&+\|\partial_x^3u_2\|_{\infty}+\|\langle x\rangle^{1/4}\partial_x^2u_2\|_{\infty}+\|\partial_x u_2\|_{\infty}(\|u_1\|_{\infty}+\|u_2\|_{\infty}\big)\Big)
\end{split}
\end{equation*}
with $M, L\in L^{\infty}([0,T])$. Therefore
\begin{equation*}
\underset{[0,T]}{\sup} \int w^2(x,t)\phi_N(x,t)\,dx \le c\Big( \int  w^2(x,0)\phi_N(x,0)\,dx +\int_0^T L(t)\,dt \Big)\,e^{\int_0^T M(t)\,dt}.
\end{equation*}
which basically yields the desired result.

\section{Proofs of Theorem \ref{theorem3} and Theorem \ref{theorem4}}

\noindent\underline{Proof of Theorem \ref{theorem3}}
\vskip5pt

To simplify the exposition and illustrate the argument of proof we restrict ourselves to consider the most difficult case 
$P(u,\partial_x u, \partial_x^2 u, \partial_x^3 u)=\partial_x^2u\partial_x^3 u$. Thus we have the equation
\begin{equation*}
\partial_t u -\partial_x^5 u+a\,\partial_x^2u\partial_x^3u=0,\hskip15pt a\in\R.
\end{equation*}

Now we follow the argument given in the proof of Theorem \ref{theorem1}. Then we need to consider the term
\begin{equation*}
I=\int \partial_x^2u\partial_x^3u u\phi_N\,dx.
\end{equation*}

By integration by parts it follows that 
\begin{equation*}
\begin{split}
I&=\frac{1}{20}\int \partial_x^5(u\,u)u\phi_N\,dx-\frac{1}{10}\int u\partial_x^5u u\phi_N\,dx-\frac12\int \partial_xu\partial_x^4u u\phi_N\,dx\\
&=I_1+I_2+I_3.
\end{split}
\end{equation*}
Hence
\begin{equation*}
I_1=-\frac{1}{20}\int \partial_x^5(u\phi_N)\,u^2\,dx
\end{equation*}
Thus
\begin{equation}\label{4.1}
|I_1|\le c\,\underset{j=0}{\overset{5}{\sum}}\|\langle x_{+}\rangle^{j/4}\partial_x^{5-j}u\|_{\infty}\int u^2\phi_N\,dx.
\end{equation}

Similarly,
\begin{equation*}
I_2\le c\|\partial_x^5 u\|_{\infty}\int u^2\phi_N\,dx,
\end{equation*}
and after integration by parts
\begin{equation*}
I_3\le \big(\|\partial_x^5u\|_{\infty}+\|\langle x_{+}\rangle^{1/4}\partial_x^4u\|_{\infty}\big)\int u^2\phi_N\,dx.
\end{equation*}

Therefore one has that
\begin{equation}\label{4.2}
|I|\le c\, \underset{j=0}{\overset{5}{\sum}}\|\langle x_{+}\rangle^{j/4}\partial_x^{5-j}u\|_{\infty}\int u^2\phi_N\,dx.
\end{equation}

Now using Theorem \ref{theorem5} and interpolation one has that for $j=0,1,\dots,6$
\begin{equation*}
\underset{[0,T]}{\sup} \|e^{(6-j)\beta x}\partial_x^j u(t)\|_2\le c,
\end{equation*}
which combined with \eqref{4.2} and Sobolev embedding yields the desired result.

\vskip3mm

\noindent\underline{Proof of Theorem \ref{theorem4}}
\vskip5pt

As in the proof of Theorem \ref{theorem3} we shall consider the most significant form of the polynomial $P(\cdot)$ in \eqref{I22b},  
$P(u,\partial_x u, \partial_x^2 u, \partial_x^3 u)=\partial_x^2u\partial_x^3 u$. Thus we consider the equation
\begin{equation*}
\partial_t u -\partial_x^5 u+a\,\partial_x^2u\partial_x^3u=0,\hskip15pt a\in\R.
\end{equation*}

Hence $w=u_1-u_2$ satisfies
\begin{equation*}
\partial_t w -\partial_x^5 w+a\,\partial_x^2u_1\partial_x^3w+a\,\partial_x^3u_2\partial_x^2w=0.
\end{equation*}

Following the argument given in the proof of Theorem \ref{theorem2} we shall estimate
\begin{equation*}
E_1=\int \partial_x^2u_1\partial_x^3w w\phi_N\,dx
\end{equation*}
and
\begin{equation*}
E_2=\int \partial_x^3u_2\partial_x^2w w\phi_N\,dx.
\end{equation*}

More precisely, we have
\begin{equation}\label{4.3}
\frac{d}{dt}\int w^2\phi_N\,dx +\epsilon \int (\partial_x^2w)^2\partial_x\phi_N\,dx\le c_0\,\int w^2\phi_n\,dx+E_1+E_2.
\end{equation}

To bound $E_2$ we use Corollary \ref{corollaryA} and \eqref{A10}
\begin{equation}\label{4.4}
\begin{split}
|E_2| &\le \tilde{c}_0 \int \partial_x^3u_2 \partial_x^2w w (1+\langle x\rangle \partial_x\phi_N)\,dx\\
&\le \tilde{c}_0 \int  \partial_x^3u_2\, \partial_x^2(u_1-u_2)(u_1-u_2)\,dx\\
&\;\;\;\;+c_{\epsilon'}\int ( \partial_x^3u_2)^2 \langle x\rangle^2 w^2\partial_x\phi_N\,dx +\epsilon' \int (\partial_x^2 w)^2\partial_x \phi_N\,dx\\
&\le M(t)+c_{\epsilon'}\int ( \partial_x^3u_2)^2 \langle x\rangle^{2+1/4} w^2 \phi_N\,dx +\epsilon' \int (\partial_x^2 w)^2\partial_x \phi_N\,dx\\
&\le M(t)+c_{\epsilon'}\|\langle x\rangle^{1+1/8}\partial_x^3u_2\|_{\infty}^2\int w^2\phi_N\,dx+\epsilon' \int (\partial_x^2 w)^2\partial_x \phi_N\,dx,
\end{split}
\end{equation}
where $0<\epsilon'\ll \epsilon$.

To control $E_1$ we write
\begin{equation*}
\begin{split}
E_1&=-\int \partial_x^3u_1\partial_x^2w w\phi_N\,dx-\int \partial_x^2u_1\partial_x^2w \partial_x w \phi_N-\int \partial_x^2 u_1 \partial_x^2 w w \partial_x\phi_N\,dx \\
&= E_1^1+E_1^2+E_1^3.
\end{split}
\end{equation*}

The bound for $E_1^1$ is similar to the one deduced above for $E_2$. For $E_1^3$ we write
\begin{equation*}
|E_1^3| \le \epsilon' \int (\partial_x^2w)^2 \partial_x \phi_N\,dx + c_{\epsilon'} \int (\partial_x u_1)^2 w^2\partial_x\phi_N\,dx,
\end{equation*}
hence a bound similar to that obtained for $|E_2|$ applies.

Finally, to estimate $E_1^2$ we write
\begin{equation*}
\begin{split}
E_1^2&=\int \partial_x^2u_1\partial_x^2 w\partial_x w\phi_N\,dx =\frac12 \int \partial_x^3u_1\partial_x w\partial_x w\phi_N\,dx
+ \frac12 \int \partial_x^2u_1\partial_x w\partial_x w\partial_x\phi_N\,dx\\
&=\frac14\int  w^2\Big[ \partial_x(\partial_x^4u_1\phi_N)+2\partial_x(\partial_x^3u_1\partial_x\phi_N)+\partial_x(\partial_x^2u_1\partial_x^2\phi_N)\Big]\,dx\\
&\;\;\;\;-\frac12\int \partial_x^3u_1\partial_x^2w w\phi_N\,dx -\frac12\int \partial_x^2u_1 w \partial_x^2 w\partial_x \phi_N\,dx\\
&= E_1^{2,1}+E_1^{2,2}+E_1^{2,3}.
\end{split}
\end{equation*}

Thus from \eqref{A11} one has that
\begin{equation*}
| E_1^{2,1}|\le \underset{j=0}{\overset{5}{\sum}}\|\langle x\rangle^{j/4} \partial_x^{5-j} u_1\|_{\infty}\int w^2\phi_N\, dx
\end{equation*}
and also 
\begin{equation*}
\begin{split}
|E_1^{2,3}|&\le \epsilon' \int (\partial_x^2 w)^2\partial_x \phi_N\, dx + c_{\epsilon'} \int  (\partial_x^2 u_1)^2 w^2 \partial_x \phi_N\, dx\\
&\le  \epsilon' \int (\partial_x^2 w)^2\partial_x \phi_N\, dx+ c_{\epsilon'} \int  (\partial_x^2 u_1)^2 \langle x\rangle^{1/4} w^2 \phi_N\, dx.
\end{split}
\end{equation*}

Finally an argument similar to that given in \eqref{4.4} shows that 
\begin{equation*}
|E_1^{2,2}| \le M(t) +c_{\epsilon'}\int (\partial_x^3 u_1)^2 \langle x\rangle^{2+1/4} w^2\phi_N\,dx+\epsilon' \int (\partial_x^2 w)^2\partial_x \phi_N\, dx.
\end{equation*}
 Inserting these estimates in \eqref{4.3} one gets the desired result.

 \section{Proof of Theorem \ref{theorem5}}

\noindent\underline{Proof of Theorem \ref{theorem5}} 
\vskip5pt

We shall follow Kato's approach in \cite{kato-83} and define for $\beta>0$
\begin{equation}\label{5.1}
\varphi_{\delta}(x)=\frac{e^{\beta x}}{1+\delta e^{\beta x}} \text{\hskip10pt for \hskip10pt} \delta\in(0,1), \;\;\;\delta\ll 1.
\end{equation}

Thus one has
\begin{equation}\label{5.2}
\varphi_{\delta}\in L^{\infty}(\R) \text{\hskip10pt and \hskip10pt} \|\varphi_{\delta}\|_{\infty}= \frac{1}{\delta}.
\end{equation}
\begin{equation}\label{5.3}
0\le \partial_x \varphi_{\delta}(x)=\frac{\beta \, e^{\beta x}}{(1+\delta e^{\beta x})^2} \le \beta \,\varphi_{\delta}(x),
\end{equation}
\begin{equation}\label{5.4}
\partial_x^2 \varphi_{\delta}(x)=\frac{\beta^2 \, e^{\beta x}(1-\delta e^{\beta x})}{(1+\delta e^{\beta x})^3} ,
\end{equation}
then
\begin{equation}\label{5.5}
|\partial_x^2\varphi_{\delta}(x)|\le \beta^2 \frac{e^{\beta x}}{(1+\delta e^{\beta x})^2}.
\end{equation}
\begin{equation}\label{5.6}
\partial_x^3 \varphi_{\delta}(x)=\frac{\beta^3 \, e^{\beta x}(1-4\delta e^{\beta x}+\delta^2e^{2\beta x})}{(1+\delta e^{\beta x})^4} ,
\end{equation}
so
\begin{equation}\label{5.7}
|\partial_x^3\varphi_{\delta}(x)|\le 2\beta^3 \frac{e^{\beta x}}{(1+\delta e^{\beta x})^2}.
\end{equation}
and
\begin{equation}\label{5.8}
|\partial_x^j\varphi_{\delta}(x)|\le c_j\beta^j \frac{e^{\beta x}}{(1+\delta e^{\beta x})^2}, \;\;\;\; j=1, 2, 3, 4, 5.
\end{equation}

Also we have that
\begin{equation}\label{5.9}
0\le \frac{(\partial_x^3\varphi_{\delta})^2}{\partial_x\varphi_{\delta}}\le 4\beta^5 \frac{e^{\beta x}}{(1+\delta e^{\beta x})^2}.
\end{equation}

Therefore
\begin{equation}\label{5.10}
\frac32|\partial_x^5\varphi_{\delta}(x)|+\frac{25}{4(5-\epsilon)} \frac{(\partial_x^3\varphi_{\delta})^2}{\partial_x\varphi_{\delta}}\le c_0\beta^5  \frac{e^{\beta x}}{(1+\delta e^{\beta x})^2}
\le c_0\beta^5 \varphi_{\delta}(x).
\end{equation}

Moreover
\begin{equation}\label{5.11}
\varphi_{\delta}(x)\le \varphi_{\delta'}(x) \quad x\in\R \text{\hskip10pt if \hskip10pt} 0<\delta'<\delta
\end{equation}
and
\begin{equation}\label{5.12}
\underset{\delta\downarrow 0}{\lim}  \, \varphi_{\delta}(x)=e^{\beta x}. 
\end{equation}

As in Theorem \ref{theorem3} and Theorem \ref{theorem4} we shall consider the most relevant case in \eqref{I22b}.
\begin{equation}\label{5.13}
P(u,\partial_xu, \partial_x^2u, \partial_x^3u)= a\partial_x^2 u\partial_x^3 u, \hskip10pt a\in\R,
\end{equation}
to get the equation
\begin{equation}\label{5.14}
\partial_tu- \partial_x^5u+ a\partial_x^2 u\partial_x^3 u=0.
\end{equation}

We employ an argument similar  to that exposed in \eqref{3.6}. Indeed, we multiply equation \eqref{5.14} by $u\varphi_{\delta}$ and integrate by parts. Then we use the
Cauchy-Schwarz and Young inequalities and the property \eqref{5.10}, to obtain the estimate
\begin{equation}\label{5.14b}
\begin{split}
&2\big(\int \partial_t u \,u\varphi_{\delta}\,dx-\int\partial_x^5\,u\varphi_{\delta}\,dx\big)\\
&=\frac{d}{dt} \int u^2\,\varphi_{\delta}\,dx+ 5\int (\partial_x^2u)^2\partial_x\varphi_{\delta}\,dx+5\int u\partial_x^2u\partial_x^3\varphi_{\delta}\,dx-\frac32 \int u^2\partial_x^5\varphi_{\delta}\,dx\\
&\ge \frac{d}{dt} \int u^2\,\varphi_{\delta}\,dx+ \epsilon \int (\partial_x^2u)^2\partial_x\varphi_{\delta}\,dx- \int u^2\,\Big(\frac32|\partial_x^5\varphi_{\delta}|+
\frac{25}{4(5-\epsilon)} \frac{( \partial_x^3\varphi_{\delta})^2}{\partial_x\varphi_{\delta}}\Big) \,dx\\
&\ge \frac{d}{dt} \int u^2\,\varphi_{\delta}\,dx+ \epsilon\int (\partial_x^2u)^2\partial_x\varphi_{\delta}\,dx-c_0 \beta^5 \int u^2\,\varphi_{\delta}\,dx
\end{split}
\end{equation}
with $\epsilon\in[0,1)$, $\,\epsilon\ll 1$ and $c_0>0$.  With this estimate we deduced that
\begin{equation}\label{5.15}
\begin{split}
\frac{d}{dt}\int u^2 \varphi_{\delta}(x)\,dx&+ \epsilon \int (\partial_x^2 u)^2 \partial_x\varphi_{\delta}(x)\,dx\\
&\le c_0\beta^5 \int u^2\varphi_{\delta}(x)\,dx +|a\int \partial_x^2 u\partial_x^3u u\varphi_{\delta}(x)\, dx|.
\end{split}
\end{equation}

Next we estimate the last term of \eqref{5.15}. We integrate by parts and write
\begin{equation*}
\begin{split}
\int \partial_x^2u\partial_x^3u u\varphi_{\delta}(x)\,dx&=\frac{1}{20}\int \partial_x^5(u^2) u \varphi_{\delta}(x)\,dx
-\frac{1}{10} \int u \partial_x^5 u u \varphi_{\delta}(x)\,dx-\frac12\int \partial_x u \partial_x^4 u u\varphi_{\delta}(x)\,dx\\
&=E_1+ E_2+ E_3.
\end{split}
\end{equation*}

Thus one has
\begin{equation*}
E_1= -\frac{1}{20}\int u^2 \partial_x^5(u \varphi_{\delta}(x))\,dx.
\end{equation*}

Therefore by \eqref{5.1}-\eqref{5.8}
\begin{equation*}
|E_1|\le c\,\underset{j=0}{\overset{5}{\sum}} \beta^j \|\partial_x^{5-j}u(t)\|_{\infty} \int u^2\varphi_{\delta}(x)\,dx.
\end{equation*}

Also
\begin{equation*}
|E_3| \le c\,(\|\partial_x^5 u\|_{\infty}+\beta \|\partial_x^4 u\|_{\infty})\,\int u^2 \varphi_{\delta}(x)\,dx
\end{equation*}
and
\begin{equation*}
|E_2| \le \|\partial_x^5 u\|_{\infty}\,\int u^2 \varphi_{\delta}(x)\,dx.
\end{equation*}

Inserting these estimates in \eqref{5.15} it follows that
\begin{equation*}
\frac{d}{dt}\int u^2 \varphi_{\delta}(x)\,dx\le c_0\big(\beta^5+\underset{j=0}{\overset{5}{\sum}} \beta^j \|\partial_x^{5-j}u(t)\|_{\infty}\big) \int u^2\varphi_{\delta}(x)\,dx
\end{equation*}
which implies that
\begin{equation}\label{5.16}
\begin{split}
\underset{[0,T]}{\sup} \int u(x,t) \varphi_{\delta}(x)\,dx &\le  \int u_0(x) \varphi_{\delta}(x)\,dx\; e^{\int_0^T N(t)\,dt}\\
&\le  \int u_0(x) \varphi_0(x) \,dx\; e^{\int_0^T N(t)\,dt}
\end{split}
\end{equation}
with
\begin{equation*}
N(t) =c_0\big(\beta^5+\underset{j=0}{\overset{5}{\sum}} \beta^j \|\partial_x^{5-j}u(t)\|_{\infty}\big).
\end{equation*}

Since the right hand side of \eqref{5.16} is independent of $\delta$ taking $\delta\downarrow 0$ we obtain the desired result.

We shall notice that in the argument above we assumed the solution sufficiently smooth to perform the integration by parts, otherwise we consider the IVP  associated to the equation \eqref{5.14} with regularized initial data as was done in the proof of Theorem \ref{theorem1}.

\vskip.2in

\section*{Acknowledgments}
P. I. was supported by DIME Universidad Na\-cio\-nal de Co\-lom\-bia-Me\-de\-ll\'in, grant 
2010\-1001\-1032. F. L. was partially supported
by CNPq and FAPERJ/Brazil. G. P. was  supported by a NSF grant  DMS-1101499. 
\vskip.2in


\end{document}